\documentclass[amssymb,amstex,amsmath,11pt,a4paper]{amsart}
\usepackage{amsmath,amssymb,latexsym,amsfonts,url}
\usepackage{graphicx}
\usepackage{epsfig}

\newtheorem{thm}{Theorem}[section]
\newtheorem*{theorem}{Theorem}

\newtheorem{lem}[thm]{Lemma}

\newtheorem{rem}[thm]{Remark}

\newcommand{\im}{{\rm Im}}

\newenvironment{pf}{\noindent {\it Proof.}}{\\}

\numberwithin{equation}{section}

\begin{document}
\title[Spacelike capillary surfaces]{Spacelike capillary surfaces \\in the Lorentz-Minkowski space}
\author[J. Pyo and K. Seo]{Juncheol Pyo and Keomkyo Seo}

\begin{abstract}
For a compact spacelike constant mean curvature surface with nonempty boundary in the three-dimensional Lorentz-Minkowski space, we introduce a rotation index of the lines of curvature at the boundary umbilic point, which was developed by Choe \cite{Choe}. Using the concept of the rotation index at the interior and boundary umbilic points and applying the Poincar\'{e}-Hopf index formula, we prove that a compact immersed spacelike disk type capillary surface with less than $4$ vertices in a domain of $\Bbb L^3$ bounded by (spacelike or timelike) totally umbilic surfaces is part of a (spacelike) plane or a hyperbolic plane. Moreover we prove that the only immersed spacelike disk type capillary surface inside de Sitter surface in $\Bbb L^3$ is part of (spacelike) plane or a hyperbolic plane. \\

\end{abstract}
\maketitle
\section{Introduction}

Spacelike surfaces with constant mean curvature (CMC) in the three-dimensional Lorentz-Minkowski space have been studied for a long time. Besides the importance of spacelike CMC surface in mathematics, such surfaces have played an important role in the study of general relativity. (See \cite{CBY} and \cite{MT} for a survey.) In \cite{ALP}, using integral formulas for compact spacelike CMC surfaces in $\Bbb L^3$, Al\'{i}as, L\'{o}pez, and Pastor proved that the only immersed compact spacelike CMC surfaces in $\Bbb L^3$ spanning a circle are the (spacelike) planar disks and the hyperbolic caps. Moreover this uniqueness result was generalized into the $n$-dimensional Lorentz-Minkowski space by Al\'{i}as and Pastor \cite{AP98}. One year later, Al\'{i}as and Pastor \cite{AP99} introduced a variational problem for spacelike surfaces in $\Bbb L^3$ whose critical points are indeed spacelike CMC surfaces intersecting a given support surface of a constant hyperbolic angle. For these spacelike CMC surfaces with free boundary in $\Bbb L^3$, they were able to prove the following.
\begin{theorem} [\cite{AP99}]
The only immersed spacelike CMC surfaces in $\Bbb L^3$ with (spacelike) planar or hyperbolic support surfaces are the planar disks ($H=0$) and the hyperbolic caps ($H \neq 0$).
\end{theorem}

On the other hand, analogous problems for CMC surfaces in the Euclidean space concerning planar disks and spherical caps have been studied as well (\cite{ALP99}, \cite{EBMR}, \cite{Koiso}, \cite{LM95}, \cite{LM96}, \cite{RR}). In particular, it is well-known that a capillary disk in a ball of the three-dimensional Euclidean space must be totally umbilic (\cite{Nitsche}, \cite{RS}). This is called a Nitsche's theorem. Here a {\it capillary surface} $M$ in a domain $U$ is a constant mean curvature surface which meets $\partial U$ in a constant contact angle along $\partial M \cap \partial U$. Physically capillary surfaces arise as the surface of an incompressible liquid in a container.(See \cite{Finn} and references therein.) In 2002, Choe \cite{Choe} showed that if a compact immersed disk type capillary surface ($H \neq0)$ in a domain bounded by planes or spheres in $\Bbb R^3$ has less than four vertices on its boundary, then the surface must be spherical. Because a regular capillary disk has no vertices on its boundary, Choe's result can be thought as a generalization of Nitsche's theorem. Motivated by this, we investigate a compact immersed spacelike capillary surface with vertices in $\Bbb L^3$ bounded by (spacelike or timelike) totally umbilic surfaces.

\begin{figure}\label{figure}
\begin{center}
\includegraphics[height=3.5cm, width=12cm]{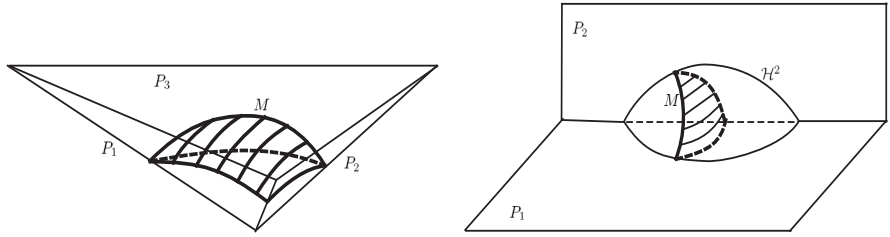}
\end{center}
\caption{Spacelike capillary surface $M$ with three vertices inside a domain bounded by three (spacelike or timelike) planes $P_{i}$, $i=1,2,3.$ (Left);
Spacelike capillary surface $M$ with three vertices inside a domain bounded by a spacelike plane $P_1$, a timelike plane $P_2$ and a hyperbolic plane $\mathcal{H}^2$ (Right).}
\end{figure}

In order to deal with spacelike CMC surfaces with vertices, we introduce a rotation index of the lines of curvature at the boundary umbilic point, which was first studied by Choe \cite{Choe}. Using the concept of the rotation index at the interior and boundary umbilic points and applying the Poincar\'{e}-Hopf index formula, we prove
\begin{theorem}
Let $M \subset \Bbb L^3$ be a compact spacelike immersed disk type CMC surface which is $C^{2,\alpha}$ up to and including $\partial M$ and whose boundary is $C^{2,\alpha}$ up to and including its vertices. Suppose that each regular component of $\partial M$ is a line of curvature. If the number of vertices of $M$ with angle $<\pi$ is less than or equal to $3$, then $M$ is part of a (spacelike) plane or a hyperbolic plane.
\end{theorem}
As a consequence of this theorem, we obtain the following uniqueness theorem. (See the Figure \ref{figure}.)
\begin{theorem}
Let $U \subset \Bbb L^3$ be a domain bounded by (spacelike or timelike) totally umbilic surfaces in $\Bbb L^3$  and let $M$ be a compact spacelike immersed disk type capillary surface in $U$ which is $C^{2,\alpha}$ up to and including $\partial M$ and whose boundary is $C^{2,\alpha}$ up to and including its vertices. If $M$ has less than $4$ vertices with angle $<\pi$, then $M$ is part of a (spacelike) plane or a hyperbolic plane.
\end{theorem}

Our theorems can be regarded as an extension of Al\'{i}as and Pastor \cite{AP99}, since they showed the uniqueness theorem for regular spacelike capillary surfaces which have no vertices. As another application, we prove that the only immersed spacelike disk type capillary surface inside a de Sitter surface in $\Bbb L^3$ is part of a (spacelike) plane ($H = 0$) or a hyperbolic plane ($H\neq 0$). (See Theorem \ref{thm:de Sitter}.)

\section{Preliminaries}
Let $\mathbb{L}^3$ be the three dimensional
Lorentz-Minkowski space, that is, the real vector space
$\mathbb{R}^3$ endowed with the Lorentz-Minkowski metric $\langle,\rangle$,
where $\langle,\rangle = d{x_1}^2 +d{x_2}^2 - d{x_3}^2$ and $x_1, x_2, x_3$
are the canonical coordinates of $\mathbb{R}^3$. If $M \subset \Bbb L^3$ is an embedded connected spacelike surface, we shall denote by $N_M$ the unique future-directed unit normal timelike vector field on $M$. Here we call a vector $v \in \Bbb L^3$ {\it future-directed} if $v$ has the same orientation as $(0,0,1)\in \Bbb L^3$.  We say that a
vector $v \in \mathbb{R}^3 \setminus \{0\}$ is {\it spacelike, timelike or
lightlike} if $|v|^2 = \langle v,v\rangle$ is positive, negative or zero,
respectively. The zero vector $0$ is spacelike by convention. A
plane in $\mathbb{L}^3$ is said to be {\it spacelike, timelike or lightlike} if the
normal vector of the plane is  timelike, spacelike, or lightlike, respectively. An immersed
 surface $\Sigma \subset \mathbb{L}^3$ is called {\it spacelike} if every tangent
 plane is a spacelike.
We now give some examples of spacelike and timelike surfaces.
\begin{itemize}
\item[(i)] The horizontal plane $\{x_{3}=c\}$ for a constant $c\in \Bbb R$ is spacelike and the vertical plane $\{a x_1+b x_2=0\}$ is timelike for any constants $a,b \in \Bbb R$ except $a=b=0$.
\item[(ii)] The hyperbolic plane $$\mathcal{H}^2(-c)=\{x=(x_1, x_2, x_3)\in \Bbb L^3: \langle x,x \rangle=-c^2, x_{3}>0\}$$ is a spacelike surface for a positive constant $c\in \Bbb R$. The unit normal vector is itself for each point on the hyperbolic plane.
\item[(iii)] The de Sitter surface is defined as $$\mathcal{S}^{2}(c)=\{x\in \Bbb L^3: \langle x,x \rangle=c^2\}$$ for a positive constant $c\in \Bbb R$. Note that the de Sitter surface is timelike and the unit normal vector is also itself for each point on the de Sitter surface.
\end{itemize}
Let $M \subset \Bbb L^3$ be a spacelike or timelike surface.
 A point $p\in M$ is called \textit{umbilic} if for any $\xi_1, \xi_2 \in T_p M$,
  $$II_{p}(\xi_{1},\xi_{2})=\lambda(p)\langle \xi_{1},\xi_{2} \rangle,$$
  that is, the second fundamental form $II$ is proportional to the first fundamental form. In case where the immersion is spacelike,
  this is equivalent to saying that two principle curvature are equal at $p$. A surface is called \textit{totally umbilic} if any point is umbilic. The (spacelike or timelike) totally umbilic surfaces in the three dimensional Lorentz-Minkowski space are classified as follows:
\begin{thm}\cite[p.116]{ONeill} \label{thm:oneill}
The only totally umbilic surfaces in $\Bbb L^3$ are planes, hyperbolic planes, and de Sitter surfaces.
\end{thm}

Throughout this paper, we shall use two different Lorentzian timelike angles in the three-dimensional Lorentz-Minkowski space in addition to a usual Lorentzian spacelike angle \cite{Ratcliffe}. If $u$ and $v$ are future-directed timelike vectors in $\Bbb L^3$, then the {\it Lorentzian timelike angle} between $u$ and $v$ is defined to be a unique nonnegative real number $\beta(u, v)$ such that
\begin{eqnarray*}
\langle u, v \rangle = |u| |v| \cosh \beta (u, v) ,
\end{eqnarray*}
where $|w|$ means the absolute value of $\langle w, w \rangle^{1/2}$ for a timelike vector $w \in \Bbb L^3$. In fact, this Lorentzian timelike angle between two timelike vectors was called a {\it hyperbolic angle} in \cite{AP99}.  If $u$ is a spacelike vector and $v$ is a future-directed timelike vector in $\Bbb L^3$, then the {\it Lorentzian timelike angle} between $u$ and $v$ is defined to be a unique nonnegative real number $\beta(u, v)$ such that
\begin{eqnarray*}
\langle u, v \rangle = |u| |v| \sinh \beta (u, v) .
\end{eqnarray*}
For simplicity, we will call $\beta(u, v)$ the {\it angle} between $u$ and $v$.

\section{Rotation index}
Let $X:M \hookrightarrow \Bbb L^3$ be a spacelike immersion. Consider an
isothermal coordinate $z=u+iv$ on $M$ taking values in a simply
connected domain $\Omega\subset M$. The metric of $\Omega$ induced
by the immersion $X$ can be written by $ds^2=\lambda^2|dz|^2$. Since every
spacelike surface is orientable, one can define a timelike normal vector field
$N$ on $M$ which satisfies $\langle N,N \rangle=-1$. Thus it is easy to see that $\{X_{u}, X_{v}, N\}$ is an orthogonal frame on $\Omega$ and
\begin{eqnarray*}
X_{uu} & = & \frac{\lambda_{u}}{\lambda}X_{u}-\frac{\lambda_{v}}{\lambda}X_{v}-eN \nonumber\\
X_{uv} & = & \frac{\lambda_{v}}{\lambda}X_{u}+\frac{\lambda_{u}}{\lambda}X_{v}-fN \nonumber\\
X_{vv} & = & -\frac{\lambda_{u}}{\lambda}X_{u}+\frac{\lambda_{v}}{\lambda}X_{v}-gN, \nonumber\\
\end{eqnarray*}
where $e=-\langle N_u, X_u\rangle = \langle N, X_{uu}\rangle$, $f=-\langle N_u, X_v\rangle = \langle N, X_{uv}\rangle=-\langle N_v, X_u\rangle$, and $g=-\langle N_v, X_v\rangle = \langle N, X_{vv}\rangle$. If we put $\Phi(z,\bar{z})=e-g-2if$, the Codazzi equation implies that Hopf function
$\Phi(z,\bar{z})$ is holomorphic with respect to the complex coordinate
$z$ if and only if the immersion $X$ is a spacelike CMC immersion. It is well-known that $\Phi dz^2$ is a holomorphic quadratic differential.

Since the principal curvature $\kappa$ and the infinitesimal principal vector $\begin{pmatrix}   du\\ dv\\  \end{pmatrix}$ satisfy $dN \begin{pmatrix}   du\\ dv\\  \end{pmatrix} = \kappa \begin{pmatrix}   du\\ dv\\  \end{pmatrix}$, one can see that the equation for the lines of curvature is given by
$$-fdu^2+(e-g)dudv+fdv^2=0 ,$$
which implies for $z=u+iv$, $$\im (\Phi dz^2)=0.$$ This is equivalent to
\begin{eqnarray*}
\arg \Phi+2\arg(dz)=m\pi~~(m ~~\text{an integer})
\end{eqnarray*}
or
\begin{eqnarray*}
\arg(dz)=\frac{m\pi}{2}-\frac{1}{2}\arg \Phi .
\end{eqnarray*}

The lines of curvature of a spacelike surface generate a smooth
line field except at umbilic points. They rotate sharply around the
umbilic points. Note that the umbilic points are isolated because such points are the zeros of the holomorphic function $\Phi$.
 The rotation index of the lines
of curvature at an umbilic point is defined as
\begin{eqnarray*}
I=\frac{1}{2\pi}\delta(\arg dz)=-\frac{1}{4\pi}\delta(\arg
\Phi),
\end{eqnarray*}
where $\delta$ denotes the variation if one goes once around an isolated umbilic point. Therefore
if the umbilic point $p$ is in the interior of the spacelike
surface, then $p$ is zero of $\Phi$ of order $n(\geq 1)$ and
$\delta(\arg \Phi)=2\pi n$. Thus at an interior point $p$ we have
\begin{eqnarray} \label{ineq:interior}
I(p)=-\frac{n}{2}\leq-\frac{1}{2}.
\end{eqnarray}

So far we have discussed the rotation index of the lines of curvature at the interior umbilic points.
From now on, we consider the case where the umbilic point $q$ is in the boundary of
the spacelike surface. In \cite{Choe}, Choe introduced the rotation index of the lines of curvature at a boundary umbilic point $q$.
Using his idea, we can estimate the rotation index. Let us give the brief idea of the definition of rotation index at the boundary umbilic point.
We may assume that a
neighborhood of $q$ is a conformal immersion of a half disk
$D_{h}=\{(u,v)\in\Bbb R^2 : u^2+v^2<1, v\geq 0\}$, $X: D_{h}\rightarrow M \subset \Bbb
L^3$ with the diameter $l$ of $D_{h}$ into the boundary of the
spacelike surface and $X(0)=q$. Since $X(l)$ is a line of
curvature of $M$, this line field can be extended smoothly to a
line field $L$ on the whole disk $D=\{(u,v)\in\Bbb R^2 : u^2+v^2<1\}$ by reflection through the
diameter. So one can define the rotation index of the lines of
curvature at the boundary umbilic point $q$ to be {\it half} the
rotation index of $L$ at $0$. This is independent of the choice of
the immersion $X$. Thus at a boundary umbilic point $q$ such that $\Phi$ has a zero of order $n$, one sees that
\begin{eqnarray*}
I(q)=\frac{1}{2}\Big[-\frac{1}{4\pi}\delta(\arg \Phi)\Big] = -\frac{n}{4} .
\end{eqnarray*}
A singular point of the boundary $\partial M$ of a spacelike surface $M \subset \Bbb L^3$ is called a {\it vertex}.
For an immersed CMC surface in $\Bbb R^3$, Choe gave an estimate for the rotation index at the boundary umbilic points and vertices. Since the rotation index is intrinsically defined, we shall make use of his results without proof.
\begin{lem}[\cite{Choe}] \label{lem:Choe01}
Let $M \subset \Bbb L^3$ be an immersed spacelike CMC surface which is $C^{2,\alpha}$ up to and including $\partial M$ and whose boundary is $C^{2,\alpha}$ up to and including its vertices. If the smooth components of $\partial M$ are lines of
curvature, then the following properties hold.
\begin{itemize}
\item[(a)] The boundary umbilic
points of $M$ are isolated.
\item[(b)] At a boundary umbilic point
which is not a vertex of $M$ the rotation index of lines of
curvature is not bigger than $-1/4$.
\item[(c)] At a vertex of $M$
with angle less than $\pi$, the rotation index is less than or
equal to $1/4$, and at a vertex with angle greater than $\pi$, the
rotation index is less than or equal to $-1/4$.
\end{itemize}
\end{lem}

\begin{lem}[\cite{Choe}]\label{lem:Choe02}
Assume that $M$ and $\partial M$ are the same as in Lemma \ref{lem:Choe01} and assume that $p$ is a vertex of $M$ with angle $\xi$. If $\xi < \pi$ and $p$ is a singularity of $\Phi$, then $p$ is a simple pole, and if $\xi > \pi$ then $p$ is a zero of $\Phi$.
\end{lem}

From the above rotation index estimate, we can prove the following uniqueness theorem for an immersed spacelike CMC surface in $\Bbb L^3$. The proof is based on \cite{Choe}.

\begin{thm} \label{thm:unique}
Let $M \subset \Bbb L^3$ be a compact immersed spacelike disk type CMC surface which is $C^{2,\alpha}$ up to and including $\partial M$ and whose boundary is $C^{2,\alpha}$ up to and including its vertices. Suppose that each regular component of $\partial M$ is a line of curvature. If the number of vertices of $M$ with angle $<\pi$ is less than or equal to $3$, then $M$ is part of a (spacelike) plane or a hyperbolic plane.
\end{thm}

\begin{pf}
The well-known Poincar\'{e}-Hopf theorem \cite{Hopf} says that the sum of the rotation indices of all singularities of a vector field is equal to the Euler characteristic of the surface. Therefore one sees that if $V$ is a line field on the domain $D$ with a finite number of singularities which is the pull-back under $x:D\rightarrow M$ of the lines of curvature on $M$, then the sum of the rotation indices of $V$ at the singularities in $\bar{D}$ is equal to $1$. So $M$ has a nonempty set $S$ of singularities. Moreover, the singularities of the lines of curvature on $M$ occur not only at the umbilic points but also at the vertices of $M$. Here the umbilic points correspond to the zeros of $\Phi$ and the vertices correspond to the poles or zeros of $\Phi$ by Lemma \ref{lem:Choe02}.

Now suppose the singular set $S$ is finite. Let $p_i, q_j, r_k$, and $s_l$ be the interior umbilic points, nonvertex boundary umbilic points, vertices with angle $>\pi$, and vertices with angle $<\pi$, respectively. Then from the inequality (\ref{ineq:interior}) and Lemma \ref{lem:Choe01} it follows that
\begin{eqnarray*}
\sum_{p=p_i, q_j, r_k, s_l} I(p) &\leq& \sum_i (-{1 \over 2}) + \sum_j (-{1 \over 4}) + \sum_k (-{1 \over 4}) + \sum_l ({1 \over 4})\\
&\leq& \sum_l ({1 \over 4}) \\
&\leq& {3 \over 4} ~, \mbox{ \hspace{1cm}   (by hypothesis)}
\end{eqnarray*}
which contradicts the fact that $\sum I(p) = 1$ from the Poincar\'{e}-Hopf theorem. Hence one obtains that $S$ is infinite and has an accumulation point $q$. Furthermore, if $\kappa_1$ and $\kappa_2$ are principal curvatures of $M$, then $S$ is the zero set of the continuous function $\kappa_1 - \kappa_2$, and hence $q\in S$ and $S$ is closed. However the points of $S$ except the vertices with angle $< \pi$ are also the zeros of $\Phi$. Since the zero set of the holomorphic function $\Phi$ is either open or finite, it follows that $S=M$, and therefore $M$ is totally umbilical. Thus one can conclude that $M$ is part of a (spacelike) plane or a hyperbolic plane.

\qed
\end{pf}

\section{Spacelike capillary surfaces}

Consider a domain $U \subset \Bbb L^3$ whose boundary $\partial U$ is an embedded connected spacelike or timelike surface piecewisely. A spacelike {\it capillary surface} $M$ in a domain $U\subset \Bbb L^3$ is an immersed spacelike CMC surface which meets $\partial U$ at a constant contact angle along $\partial M$. If $\partial U$ is a piecewise smooth surface, then we may assume the constant angles to be distinct on each smooth component of $\partial U$.

Let $M \subset \Bbb L^3$ be a spacelike capillary surface which meets $\partial U$ at a constant contact angle $\beta$. We shall denote by $\tau$ the positively oriented unit tangent vector field along $\partial M$ and denote by $N$ the timelike unit normal vector field on $M$. Clearly $\{\tau, N, \nu= -\tau \wedge N\}$ is trihedra along $\partial M$. Here $u \wedge v$ denotes the vector product of two vectors $u, v \in \Bbb L^3$ which is defined to be the unique vector $u \wedge v \in \Bbb L^3$ such that
\begin{eqnarray*}
\langle u \wedge v, w \rangle = \det (u, v, w)
\end{eqnarray*}
for any $w \in \Bbb L^3$ \cite{ALP}. Clearly $\nu = -\tau \wedge N$ is the inward pointing unit conormal vector field along $\partial M$. Choose a regular piece $\Sigma$ of $\partial U$. Then $\Sigma$ is an embedded connected spacelike or timelike surface. As before, one may construct trihedra $\{\tau, N_\Sigma, \nu_\Sigma \}$ along $\partial M$, where $N_\Sigma$ is the unit normal vector field on $\Sigma$ and $\nu_\Sigma$ is the inward pointing unit conormal vector field along $\partial \Sigma$ which is given by $\nu = -\tau \wedge N_\Sigma$. For these two trihedra  $\{\tau, N, \nu\}$ and $\{\tau, N_\Sigma, \nu_\Sigma \}$, we have the following equations:\\
(i) If $\Sigma$ is a spacelike surface,
\begin{eqnarray} \label{trihedra:spacelike}
\left\{ \begin{array}{ll} \nu=\cosh\beta\nu_{\Sigma}+\sinh\beta N_{\Sigma}, \\ N=\sinh\beta\nu_{\Sigma}+\cosh\beta N_{\Sigma}.\end{array}\right.
\end{eqnarray}
(ii) If $\Sigma$ is a timelike surface,
\begin{eqnarray} \label{trihedra:timelike}
\left\{ \begin{array}{ll} \nu=\sinh\beta\nu_{\Sigma}+\cosh\beta N_{\Sigma}, \\ N=\cosh\beta\nu_{\Sigma}+\sinh\beta N_{\Sigma}.\end{array}\right.
\end{eqnarray}
In case where $\Sigma$ is a spacelike surface, Al\'{i}as and Pastor \cite{AP99} also used the equation (\ref{trihedra:spacelike}) in which the constant contact angle $\beta$ only differs by a minus sign from ours.

When the ambient space is a Euclidean space, the following Terquem-Joachimsthal theorem is well known.
\begin{theorem}[\cite{Spivak}]
Let $c$ be a curve in $M_1 \cap M_2 \subset \Bbb R^3$ which is a line of curvature in $M_1$. Then $c$ is a line of curvature in $M_2$ if and only if $M_1$ and $M_2$ intersect at a constant contact angle along $c$.
\end{theorem}

This theorem can be generalized into the three-dimensional Lorentz-Minkowski space as follows. It should be mentioned that Al\'{i}as and Pastor \cite{AP99} proved this lemma for a spacelike surface $\Sigma$. The method we use here is a modification of \cite{AP99}. For the sake of completeness we give the proof.
\begin{lem} \label{lem:Joachimsthal}
Let $M \subset \Bbb L^3$ be an immersed spacelike CMC surface and let $\Sigma \subset \Bbb L^3$ be a (spacelike or timelike) totally umbilic surface. Suppose that $M$ meets $\Sigma$ at a constant contact angle along $\partial M \cap \Sigma$. Then each smooth component of $\partial M\cap \Sigma$ is a line of curvature of $M$.
\end{lem}
\begin{pf}
Choose a point $p$ on a smooth component of $\partial M \cap \Sigma$. It suffices to show that the intersection of a local neighborhood of $p$ with $\partial M \cap \Sigma$ is a line of curvature of $M$. Let $X:D_h \rightarrow M \subset \Bbb L^3$ be a conformal immersion of a half disk $D_h = \{(u, v) \in \Bbb R^2 : u^2+v^2 < 1, v \geq 0\}$ into $M$, which maps the diameter $l$ of $D_h$ into $\partial M$ and $X(0)=p$. Let $z=u+iv$ be the usual coordinate in $\Bbb C$. Then the metric on $M$ is written by $ds^2=\lambda^2|dz|^2$ for a positive smooth function $\lambda = \lambda (z)$. One can write the unit tangent vector field $\tau$ and the inward-pointing unit conormal vector field $\nu$ along the smooth boundary containing $p$ are given by $\tau=\lambda^{-1}\partial_{u}$ and $\nu=\lambda^{-1}\partial_{v}$. By Theorem \ref{thm:oneill}, we have four possible cases for $\Sigma$: a spacelike plane, a hyperbolic plane, a timelike plane, and a de Sitter surface.

When $\Sigma$ is a spacelike plane, the normal vector field of $\Sigma$ is given by $N_{\Sigma}=\overrightarrow{a}$.
\begin{eqnarray*}
II(\tau,\nu) & = & -\langle\bar{\nabla}_{\tau}N, \nu \rangle=\langle N, \bar{\nabla}_{\tau}\nu \rangle
\\
&=& \cosh\beta\langle\bar{\nabla}_{\tau}\nu_{\Sigma}, N \rangle +\sinh\beta\langle\bar{\nabla}_{\tau}\overrightarrow{a}, N \rangle \\
&=& \cosh\beta\sinh\beta\langle\bar{\nabla}_{\tau}\nu_{\Sigma}, \nu_{\Sigma} \rangle +{\cosh^2\beta}\langle\bar{\nabla}_{\tau}\nu_{\Sigma}, \overrightarrow{a} \rangle \\
&=& \frac{1}{2}\cosh\beta\sinh\beta\tau\langle\nu_{\Sigma}, \nu_{\Sigma} \rangle -{\cosh^2\beta}\langle\nu_{\Sigma}, \bar{\nabla}_{\tau}\overrightarrow{a} \rangle \\
&=& 0,
\end{eqnarray*}
where $\bar{\nabla}$ denotes the connection of $\Bbb L^3$.

When $\Sigma$ is a de Sitter surface, the normal vector field of $\Sigma$ is given by $N_{\Sigma}=X$.
\begin{eqnarray*}
II(\tau,\nu) & = & -\langle\bar{\nabla}_{\tau}N, \nu \rangle=\langle N, \bar{\nabla}_{\tau}\nu \rangle
\nonumber\\
&=& \sinh\beta\langle\bar{\nabla}_{\tau}\nu_{\Sigma}, N \rangle +\cosh\beta\langle\bar{\nabla}_{\tau}X, N \rangle \\
&=& \cosh\beta\sinh\beta\langle\bar{\nabla}_{\tau}\nu_{\Sigma}, \nu_{\Sigma} \rangle +{\sinh^2\beta}\langle\bar{\nabla}_{\tau}\nu_{\Sigma}, X \rangle \\
&=& \frac{1}{2}\cosh\beta\sinh\beta\tau\langle\nu_{\Sigma}, \nu_{\Sigma} \rangle -{\sinh^2\beta}\langle\nu_{\Sigma}, \tau\rangle\\
&=& 0
\end{eqnarray*}
When $\Sigma$ is a hyperbolic or timelike plane, the proof is similar to the case where $\Sigma$ is a spacelike plane or a de Sitter surface. \qed
\end{pf}

\begin{thm} \label{thm:cap}
Let $U \subset \Bbb L^3$ be a domain bounded by (spacelike or timelike) totally umbilic surfaces in $\Bbb L^3$ and let $M$ be a compact immersed spacelike disk type capillary surface in $U$ which is $C^{2,\alpha}$ up to and including $\partial M$ and whose boundary is $C^{2,\alpha}$ up to and including its vertices. If $M$ has less than $4$ vertices with angle $<\pi$, then $M$ is part of a (spacelike) plane or a hyperbolic plane.
\end{thm}

\begin{pf}
From Lemma \ref{lem:Joachimsthal}, we obtain that each smooth component of $\partial M$ is a line of curvature of $M$. Hence the conclusion follows from Theorem \ref{thm:unique}.
\qed
\end{pf}

\begin{rem}
In case where a smooth component of $\partial U$ is a lightlike plane, one cannot expect a similar equation as  (\ref{trihedra:spacelike}) or (\ref{trihedra:timelike}) since $N$ is lightlike vector. Thus the proof of Theorem \ref{thm:cap} does not work in this case.
\end{rem}
\begin{rem}
 The number of vertices of $M$ with angle $<\pi$ in
    Theorem \ref{thm:cap} is sharp. Let $\mathcal{C}$ be the
    Lorentzian catenoid which is a spacelike surface of revolution.(\cite{Kobayashi}, \cite{LLS})
    Consider $M\subset \mathcal{C}$ be a compact part bounded
    by two parallel horizontal spacelike planes which are perpendicular
    to the axis of  $\mathcal{C}$ and two vertical timelike
    planes containing the axis of $\mathcal{C}$ with angle
    $\theta\in (0,\pi)$. Then $M$ is a compact embedded disk
    type spacelike CMC $(H=0)$ capillary surface with four
    vertices with angle $<\pi$ at each of which the rotation
    index equals ${1}/{4}$, which is neither part of a (spacelike) plane nor a hyperbolic plane.

\end{rem}

As an immediate consequence of Theorem \ref{thm:cap}, one can obtain the following theorem which is a generalization of \cite{AP99}.

\begin{thm} \label{thm:de Sitter}
The only spacelike immersed disk type capillary surface inside a de Sitter surface in $\Bbb L^3$ is a planar disk ($H = 0$) or a hyperbolic disk ($H\neq 0$).
\end{thm}


\vspace{1cm}

\noindent Juncheol Pyo\\
School of Mathematics, Korea Institute for Advanced Study, 207-43, Cheongnyangni 2-dong, Dongdaemun-gu, Seoul, 130-722, Korea\\
{\tt E-mail:jcpyo@kias.re.kr}\\

\bigskip
\noindent Keomkyo Seo\\
Department of Mathematics, Sookmyung Women's University,
Hyochangwongil 52, Yongsan-ku, Seoul, 140-742, Korea\\
{\tt E-mail:kseo@sookmyung.ac.kr}\\


\begin{thebibliography}{123}

\bibitem{ALP99} L. Al\'{i}as, R. L\'{o}pez, B. Palmer, {\em Stable constant mean curvature surfaces with circular boundary}, Proc. Amer. Math. Soc. {\bf 127} (1999), no. 4, 1195-1200.
\bibitem{ALP} L. Al\'{i}as, R. L\'{o}pez, J. Pastor, {\em Compact spacelike surfaces with constant mean curvature in the Lorentz-Minkowski 3-space}, Tohoku Math. J. {\bf 50} (1998), 491-501.
\bibitem{AP98} L. Al\'{i}as, J. Pastor, {\em Constant mean curvature spacelike hypersurfaces with spherical boundary in the Lorentz-Minkowski space}, J. Geom. Phys. {\bf 28} (1998), no. 1-2, 85-93.
\bibitem{AP99} L. Al\'{i}as, J. Pastor, {\em Spacelike surfaces of constant mean curvature with free boundary in the Minkowski space}, Class. Quantum Grav. {\bf 16} (1999), 1323-1331.
\bibitem{Choe} J. Choe, {\em Sufficient conditions for constant mean curvature surfaces to be round}, Math. Ann. {\bf 323} (2002), no. 1, 143-156.
\bibitem{CBY} Y. Choquet-Bruhat, J. York, {\em The Cauchy problem. General relativity and gravitation, Vol. 1}, Plenum, New York-London, 1980, 99-172.
\bibitem{EBMR} R. Earp, F. Brito, W. Meeks, H. Rosenberg, {\em Structure theorems for constant mean curvature surfaces bounded by a planar curve}, Indiana Univ. Math. J. {\bf 40} (1991), no. 1, 333-343.
\bibitem{Finn} R. Finn, {\em Equilibrium capillary surfaces}, Grundlehren der Mathematischen Wissenschaften, {\bf 284}, Springer-Verlag, New York, 1986.
\bibitem{Hopf} H. Hopf, {\em Differential geometry in the large}, Lecture Notes in Mathematics {\bf 1000}, Springer-Verlag, Berlin, 1989.
\bibitem{Kobayashi} O. Kobayashi, {\em Maximal surfaces in the $3$-dimensional Minkowski space $L^{3}$},  Tokyo J. Math. {\bf 6} (1983), no. 2, 297-309.
\bibitem{Koiso} M. Koiso, {\em Symmetry of hypersurfaces of constant mean curvature with symmetric boundary}, Math. Z. {\bf 191} (1986), no. 4, 567-574.
\bibitem{LLS} F. L\'{o}pez, R. L\'{o}pez, R. Souam, {\em Maximal surfaces of Riemann type in Lorentz-Minkowski space $\Bbb L^3$}, Michigan Math. J. {\bf 47} (2000), no. 3, 469-497.
\bibitem{LM95} R. L\'{o}pez, S. Montiel, {\em Constant mean curvature discs with bounded area}, Proc. Amer. Math. Soc. {\bf 123} (1995), no. 5, 1555-1558.
\bibitem{LM96} R. L\'{o}pez, S. Montiel, {\em Constant mean curvature surfaces with planar boundary}, Duke Math. J. {\bf 85} (1996), no. 3, 583-604.
\bibitem{MT} J. Marsden, F. Tipler, {\em Maximal hypersurfaces and foliations of constant mean curvature in general relativity}, Phys. Rep. {\bf 66} (1980), no. 3, 109-139.
\bibitem{Nitsche} J. Nitsche, {\em Stationary partitioning of convex bodies}, Arch. Rational Mech. Anal. {\bf 89} (1985), no. 1, 1-19.
\bibitem{ONeill} B. O'Neill, {\em Semi-Riemannian Geometry with Application to Relativity}, Pure Appl. Math. {\bf 130}, Academic Press, 1983.
\bibitem{Ratcliffe} J. Ratcliffe, {\em Foundations of hyperbolic manifolds, Second edition}, Graduate Texts in Mathematics {\bf 149}, Springer, New York, 2006.
\bibitem{RR} A. Ros, H. Rosenberg, {\em Constant mean curvature surfaces in a half-space of $R^3$ with boundary in the boundary of the half-space}, J. Differential Geom. {\bf 44} (1996), no. 4, 807-817.
\bibitem{RS} A. Ros, R. Souam, {\em On stability of capillary surfaces in a ball}, Pacific J. Math. {\bf 178} (1997), no. 2, 345-361.
\bibitem{Spivak} M. Spivak, {\em A comprehensive introduction to differential geometry, Vol. III}, Publish or Perish, 1979.
\end{thebibliography}
\end{document}